\documentclass[11pt]{amsart}

\usepackage{mathpazo}
\usepackage[top=25mm, bottom=25mm, left=30mm, right=30mm]{geometry} 

\usepackage{amsmath,amssymb}
\usepackage{graphicx}
\usepackage{import}
\usepackage{amscd}
\usepackage{wrapfig}
\usepackage{epsfig}
\numberwithin{equation}{section}
\usepackage{tikz}

\newtheorem{theorem}{Theorem}[section]

\newcommand{\C}{{\mathbb C}}

\newcommand{\R}{{\mathbb R}}
\newcommand{\Z}{{\mathbb Z}}

\renewcommand{\H}{{\mathbb H}}

\newcommand{\teichmuller}{Teichm{\"u}ller{ }}

\makeatletter
 \let\c@theorem=\c@subsection
 \let\c@conjecture=\c@subsection
 \let\c@lemma=\c@subsection
 \let\c@proposition=\c@subsection
 \let\c@claim=\c@subsection
 \let\c@question=\c@subsection
 \let\c@criterion=\c@subsection
 \let\c@vfconj=\c@subsection
 \let\c@definition=\c@subsection
 \let\c@notation=\c@subsection
 \let\c@remark=\c@subsection
 \let\c@example=\c@subsection
 \let\c@equation=\c@subsection
 \let\c@figure=\c@subsection
 \let\c@wrapfigure=\c@subsection

\makeatother

\begin{document}

\title{Trimmed sums of twists and the area Siegel-Veech constant.}
\author{Vaibhav Gadre}
\address{School of Mathematics and Statistics, University of Glasgow, University Place, Glasgow G12 8SQ}
\email{Vaibhav.Gadre@glasgow.ac.uk}
 

\begin{abstract}
We relate trimmed sums of twists in cylinders along a typical \teichmuller geodesic to the area Siegel-Veech constant.
\end{abstract}

\maketitle

\pagestyle{empty}

\section{Introduction}
The strong law in \cite{Gad} relates trimmed sums of excursions of a random geodesic in the thin parts of moduli spaces of quadratic differentials to Siegel-Veech constants. 
This strong law is a generalisation of the Diamond-Vaaler strong law \cite{Dia-Vaa} for continued fraction coefficients: for almost every $r \in [0,1]$ the continued fraction coefficients of $r$ satisfy
\[
\lim\limits_{n \to \infty} \frac{a_1 + a_2 + \cdots + a_n - \max\limits_{j \leqslant n} a_k }{n \log n} = \frac{1}{\log 2}. 
\]

By work of Series \cite{Ser}, continued fraction coefficients can be interpreted in terms of hyperbolic geodesic rays on the modular surface $X = \H/ SL(2,\Z)$.
In the upper half-space model, consider the (vertical) hyperbolic geodesic $\gamma$ converging to $r \in [0,1] \subset \R \cup \infty = \partial \H$. 
Since $r$ is irrational, the geodesic $\gamma$ passes through infinitely many horoballs in the Ford packing.
As $\gamma$ enters and leaves the $k$-th horoball $H_k$, the distance along $\partial H_k$ between the entry and exit points is, up to a uniform additive constant, the same as the coefficient $a_k$ in the continued fraction expansion of $r$. 

The generalisation then proceeds by two steps.
Let $\Gamma$ be a non-uniform lattice in $SL(2,\R)$, such as $SL(2,\Z)$. 
First, one proves a continuous time strong law for trimmed sums of excursions of hyperbolic geodesics in cusp neighbourhoods of $X = \H/ \Gamma$.
The limit is $(2/\pi)$ times the relative volume in $X$ of the cusp neighbourhoods. 
For $SL(2,\Z)$, one can then invoke a well known asymptotic for the number $n$ of coefficients as a function of the time parameter $T$ along a random vertical geodesic on $X = \H/ SL(2,\Z)$. 
This allows us to pass from the continuous time strong law to the discrete version of Diamond-Vaaler. 
The proof of the continuous time strong law uses the exponential decay of correlations by Moore, Ratner \cite{Moo}, \cite{Rat} for the geodesic flow on $T^1 X$. 
It combines the decay of co-relations with estimates coming from the cusp geometry of $X$. 

From this point, the strong law generalises to cusp excursions of \teichmuller geodesics in the thin parts of $SL(2,\R)$-orbit closures. 
The exponential decay of correlations for the flow by Avila-Gou\"{e}zel-Yoccoz \cite{Avi-Gou-Yoc} and Avila-Resende \cite{Avi-Res} has been shown to generalise to orbit closures by Avila-Gou\"{e}zel \cite{Avi-Gou}. 
The cusp geometry is significantly more involved because the thin parts can intersect in complicated ways. 
However, regularity of the (absolutely continuous) invariant measure proved by Dozier \cite{Doz} and a counting result for short saddle connections by Eskin-Masur \cite{Esk-Mas} imply that the higher rank contributions to the estimates are asymptotically negligible. 
Analogous to the Fuchsian situation, the limit in the strong law is a fixed multiple of the Siegel-Veech constant associated to the thin part.

In this note, we state a particular version of a trimmed sum strong law not set out explicitly in \cite{Gad}. 
We set up some  some preliminaries before stating the theorem.
A quadratic differential $q$ on an oriented surface $S$ of finite type is equivalent to a half-translation structure on $S$, that is, contour integration of a square-root of $q$ defines charts from $S$ to $\C$ and the transition functions are half-translations with the form $z \to \pm z + c$. 
The charts thus define a singular flat metric on $S$. 

A cylinder for a quadratic differential $q$ is an embedded cylinder in $S$ that is a union of freely homotopic closed geodesic trajectories in its singular flat metric. 
Any cylinder can be enlarged to be maximal, that is, enlarged so that both boundary components of the cylinder contain a singularity of the flat metric.
The isotopy class of the simple closed curve given by any of the closed geodesics that sweep out a maximal cylinder is called the core curve. 

We prove the following theorem.

\begin{theorem}\label{main}
Let $S$ be an oriented surface of finite type. Let $\mu$ be a regular $SL(2,\R)$-invariant measure in the Lebesgue class on an affine invariant manifold $\mathcal{N}$ in a stratum $\mathcal{Q}(\alpha)$ of quadratic differentials on $S$. For $\mu$-almost every $q \in \mathcal{N}$, let $\phi_t(q)$ be the \teichmuller geodesic ray given by $q$. Suppose that till time $T > 0$ the ray $\phi_t$ has excursions in the thin parts of cylinders $\{\text{Cyl}_1, \ldots, \text{Cyl}_{N(T)}\}$, ordered by time.  Let $\alpha_j$ be the core curve of $\text{Cyl}_j$. Then 
\[
\lim\limits_{T \to \infty} \frac{\sum\limits_{j <  N(T)} \text{tw}(\alpha_j) - \max\limits_{j < N(T)} \text{tw}(\alpha_j) }{T \log T} = 4 c_{\text{area}}(\mathcal{N}),
\]
where $c_{\text{area}}$ is the area Siegel-Veech constant for $\mathcal{N}$ and where $\text{tw}(\alpha_k)$ is the number of twists in $\alpha$ that $\phi_t(q)$ has during its excursion in the thin part of $\text{Cyl}_j$. 
\end{theorem} 

The \teichmuller flow can be coded symbolically in a number of ways. Theorem \ref{main} thus sets up the intriguing prospect of a direct computer verification/computation of area Siegel-Veech constants. These constants are important for their relationship with other dynamical quantities such as the Lyapunov spectrum for the \teichmuller flow \cite[Theorem 1]{Esk-Kon-Zor}.  

\subsection{Acknowledgements:} I thank V. Delecroix and S. Schleimer for their suggestion to formulate Theorem \ref{main} relating trimmed sums of cylinder twists to area Siegel-Veech constants. 

\section{Preliminaries} 

\subsection{Quadratic differentials:} 
Let $S$ be an oriented surface of finite type. 
The \teichmuller space $\mathcal{T}(S)$ is the space of marked conformal structures on $S$. 
The mapping class group $\text{Mod}(S)$ is the group of orientation preserving diffeomorphisms of $S$. 
It acts on $\mathcal{T}(S)$ by changing the marking. 
The quotient $\mathcal{M}(S) = \mathcal{T}(S) / \text{Mod}(S)$ is the moduli space of Riemann surfaces $X$ of type $S$. 

For a Riemann surface $X$, let $\mathcal{Q}(X)$ be the set of meromorphic quadratic differentials on $X$ with simple poles at the punctures. 
If $(k_1, k_2, \dots, k_r)$ are the multiplicities of the zeros then $k_1 + k_2 + \dots + k_r = 2g-2+n$, where $n$ is the number of punctures. 
A quadratic differential is equivalent to a half-translation structure on $S$, i.e. it defines charts from $S$ to $\C$ with transition functions of the form $z \to \pm z + c$. 
The resulting flat metric has a cone singularity with angle $(k+2)\pi$ at a $k$-order zero (or with $k=-1$ for a simple pole). 
A quadratic differential is unit area if the corresponding singular flat metric has area 1. 
The space $\mathcal{Q}$ is stratified by the multiplicity of its zeros: we denote the stratum with multiplicities $\alpha = (k_1, k_2, \dots, k_r)$ by $\mathcal{Q}(\alpha)$. 
Each stratum is $\textup{Mod}(S)$ invariant. 
We will continue to denote the moduli space by $\mathcal{Q}(\alpha)$. 

The action of $SL(2,\R)$ on $\C = \R^2$ preserves the form of the transitions, that is, a half-translation is taken to a half-translation.
Hence, it descends to an action of $SL(2,\R)$ on  $\mathcal{Q}(\alpha)$.
The diagonal action gives the \teichmuller flow on $\mathcal{Q}(\alpha)$.

The compact part $SO(2,\R)$ leaves the conformal structure unchanged.
As a result, one gets isometric embeddings of $\H = SL(2,\R) / SO(2,\R)$ in $\mathcal{T}(S)$.
These are called \teichmuller discs.
We denote the \teichmuller disc determined by a quadratic differential $q$ as $\H(q)$. 

\subsection{Orbit closures:} 
Let $\alpha$ be a homology class in the first homology of $S$ relative to the singularities.
The period/ holonomy for $\alpha$ is the complex number given by integrating a square root of $q$ along a contour representing $\alpha$.
For a quadratic differential that is not a square of an abelian differential, the holonomy is only defined up to sign. 
The periods/holonomies of a fixed basis for the homology of $S$ relative to the singularities, give local co-ordinates on $\mathcal{Q}(\alpha)$. 
Eskin-Mirzakhani-Mohammadi \cite{Esk-Mir-Moh} showed that $SL(2,\R)$-orbit closures are affine sub-manifolds in the period co-ordinates. 
Eskin-Mirzakhani \cite{Esk-Mir} showed that ergodic $SL(2,\R)$-invariant measures are Lebesgue measures supported on such affine submanifolds.  

\subsection{Siegel-Veech formula:} 
Let $V$ be a $SL(2,\R)$-invariant loci of vectors (with multiplicities) in $\R^2 \setminus \{(0,0)\}$, that is, for each $q \in \mathcal{N}$, the subset $V(q)$ is a discrete subset of (weighted) non-zero vectors in $\R^2$ and the assignment $q \to V(q)$ is $SL(2,\R)$-equivariant. 
The Siegel-Veech transform given by an $SL(2,\R)$-invariant locus $V$ takes compactly supported functions on $\R^2$ to a function on $\mathcal{N}$.
For a compactly supported function $f$ on $\R^2$, it is defined as 
\[
\widehat{f}(q) = \sum\limits_{v \in V(q)} f(v) 
\]
Under some additional assumptions on the locus $V$ (see \cite[Section 2]{Esk-Mas}) one gets the Siegel-Veech formula
\[
\int\limits_{\mathcal{N}} \widehat{f} \, d\mu = c(V, \mu) \int\limits_{\R^2} f \,d \, \text{Leb} 
\] 
where $c(V,\mu)$ is a non-negative constant that is zero if and only $\widehat{f}$ is zero $\mu$-almost everywhere for all compactly supported functions $f$. 
See \cite[Section 2]{Esk-Mas} for the precise details. 

\subsection{Thin parts:} 
Let $V$ be an $SL(2,\R)$-invariant locus.
Let $\epsilon > 0$ be a fixed constant. 
The $\epsilon$-thin part of $\mathcal{N}$ given by $V$ is the set of all $q \in \mathcal{N}$ such that the subset $V(q)$ contains a (possibly weighted) vector $v$ satisfying $\Vert v \Vert^2 \leqslant \epsilon$.
To get a sensible definition, we have to assume that $\epsilon > 0$ is small enough depending only on $\mathcal{N}$ and $V$. 
With this assumption, we will suppress $\epsilon$ in the notation hereafter. 
So we will call the $\epsilon$-thin part of $\mathcal{N}$ given by $V$ as simply the $V$-thin part of $\mathcal{N}$,
We will denote it by $\mathcal{N}_{V, \epsilon}$. 

By applying the Siegel-Veech formula to the characteristic function of a ball of radius $\epsilon/ R$ about the origin in $\R^2$ and then letting $R \to \infty$ one derives the asymptotic
\[
\lim_{R \to \infty} \frac{\mu \left( \mathcal{N}_{V, \epsilon/R } \right) }{\pi \epsilon/R } = c(V, \mu).
\]
The details can be found in \cite[Section 2.1]{Esk-Mas}.

\section{The area Siegel-Veech constant and excursions in cylinder thin parts} 

\subsection{Area Siegel-Veech constants:} 
A saddle connection on a half-translation surface $(S,q)$ is a straight line segment in the flat metric embedded on its interior such that a singularity/ pole of $q$ is a point of the saddle connection if and only if it is its endpoint.
A metric cylinder $\text{Cyl}$ in $q$ is an embedded cylinder that is a union of freely homotopic closed trajectories of $q$. For our purposes, we will always consider maximal cylinders, that is, the boundary components are a concatenation of saddle connections. 
We will drop the adjective maximal from now as it is implicitly assumed.

Associated to cylinders, there are natural $SL(2,\R)$-invariant loci that use holonomies (periods) of core curves of cylinders in $q$. 
As an example,  the holonomy of a cylinder $\text{Cyl}_\alpha$ with core curve $\alpha$ can be weighted by $\text{Area}(\text{Cyl}_\alpha)$.
Since the $SL(2,\R)$ action preserves area we get an $SL(2,\R)$-invariant locus which we denote by $V_{\text{area}}$. 
By the work of Veech \cite{Vee} and Vorobets \cite{Vor}, the locus $V_{\text{area}}$ satisfies the assumptions required for the derivation of the Siegel-Veech formula.
The Siegel-Veech constant $c_{\text{area}} (\mathcal{N})$ in this case is called the area Siegel-Veech constant. 
See \cite[Section 1.6]{Esk-Kon-Zor} for more details. 

Let $\text{Cyl}_\alpha$ be a cylinder in the half-translation surface $(S,q)$ with core curve $\alpha$.
The flat length $\ell_q(\alpha)$ of $\alpha$ is exactly $\Vert v(\alpha) \Vert$, where $v(\alpha)$ is the holonomy of $\alpha$. We fix an $\epsilon > 0$ small enough so that there is an open set in $\mathcal{N}$ such that for any $q$ in this open set $\ell_q(\alpha)^2 > \epsilon$ for every cylinder $\text{Cyl}_\alpha$ in $(S,q)$. 
The cylinder thin part of $\mathcal{N}$ is then the subset of those $q$ such that $\ell_q (\alpha)^2 \leqslant \epsilon$ for a core curve $\alpha$ of some cylinder.  
The intersection of the cylinder thin part of $\mathcal{N}$ with a typical \teichmuller disc $\H(q)$ is a collection of horoballs. Given a cylinder $\text{Cyl}_\alpha$, the subset of $\H(q)$ where $\alpha$ is short is a horoball $H_\alpha$ whose point at infinity (in $S^1 = \partial \H(q)$) is given by the direction $\theta$ in which the cylinder given by $\alpha$ is vertical.

\subsection{Excursions:} 
Let $H$ be a horoball in a \teichmuller disc $\H(q)$.
Let $\pi: H \to \partial H$ be the closest point projection.
For a finite geodesic segment $\gamma$ that enters $H$ we define the excursion $E(\gamma, H)$ as 
\[
E(\gamma, H) = \ell_{\partial H} (\pi_H(\gamma))
\]
where the length on the right is measured by using the path metric on the boundary $\partial H$ of the horoball $H$. 

When a horoball $H$ is part of a collection arising from an $SL(2,\R)$-invariant locus $V$, then there might be a weight $A \neq 1$ assigned by $V$.
For example, a cylinder horoball could be weighted by the flat area of the cylinder.
In such cases, we let the excursion $E_V(\gamma, H)$ be $A$ times the length along $\partial H$ of $\pi(\gamma \cap H)$. 
In other words, $E_V(\gamma, H) = A \, E(\gamma, H)$. 

For a geodesic ray $\gamma$ let $N(T)$ be the number of cylinder horoballs encountered by $\gamma$ till time $T$. Let us number these horoballs as $H_1, H_2, \cdots, H_{N(T)}$ in the order of increasing time. Let $E_V(\gamma, T)$ be the sum 
\[
E_V(\gamma, T) = \sum\limits_{j \leqslant N(T)} E_V(\gamma, H_j). 
\] 
The strong law in \cite{Gad} states that for $\mu$-almost every $q$ in $\mathcal{N}$ we have
\[
\lim_{T \to \infty} \frac{E_V (\gamma, T) - \max\limits_{j \leqslant N(T)} E_V (\gamma,H_j) }{T \log T } = 2 \epsilon c(V, \mu).  
\]
The original statement of the strong law in \cite{Gad} was conditioned on $\mu$ being a \emph{regular} measure, that is, $\mu$ is required to satisfy a particular quantitative version of independence for multiple (non-parallel) saddle connections being short.
See \cite[Section 4.3]{Gad} for the precise description. 
In subsequent work \cite{Doz}, Dozier proved that all ergodic $SL(2,\R)$-invariant measures in the Lebesgue class of their supporting affine invariant manifolds, are regular. 
Hence, the regularity condition can now be dropped from the statement of the strong law.
In the special case that $V = V_{\text{area}}$ the area Siegel-Veech constant shows up in the limit.

Suppose that the weights assigned by $V$ are bounded above, as for example in $V_{\text{area}}$.
We will now show that the final excursion $E_V(\gamma, H_{N(T)})$ can be excluded from the trimmed sum without changing the limit.
By \cite[Corollary 3.5]{Gad}, for any constant $c> 1/2$ if $E_V(\gamma, H_{N(T)}) > T (\log T)^c$ then there is no other $k < N(T)$ such that $E_V(\gamma, H_k) > T (\log T)^c$. 
Thus $E_V(\gamma, H_{N(T)}) $ gets trimmed if it is larger than $T (\log T)^c$.
Otherwise $E_V(\gamma, H_{N(T)})$ is bounded above by $T (\log T)^c$ and hence $E_V(\gamma, H_{N(T)}) / (T \log T) \to 0$ as $T \to \infty$. 
Thus we may exclude the final excursion to conclude that for $\mu$-almost every $q$ in $\mathcal{N}$ 
\begin{equation}\label{t:strong}
\lim_{T \to \infty} \frac{ \sum\limits_{j < N(T)} E_{V_{\text{area}}} (\gamma, H_j)  - \max\limits_{j < N(T)} E_{V_{\text{area}}}(\gamma,H_j) }{T \log T } = 2 \epsilon c_{\text{area}}(\mathcal{N}).
\end{equation}
For times $T$ such that $\gamma_T$ is inside some horoball $H$, the excursion $E_V(\gamma, H_{N(T)})$ is only partial. 
The above version of the strong law allows us to restrict to trimmed sums of complete excursions excluding the final partial excursion if any.

\subsection{Twists:} We now recall from \cite[Section 2.5]{Gad-Mah-Tio} estimates relating excursions to twists in the core curves of the cylinders. 
As pointed out in \cite{Gad-Mah-Tio}, the derivation of this goes back to the work of Rafi \cite{Raf1, Raf2}. Here, we are interested in tracking the relationship of excursion to twists more carefully than it is done in \cite{Gad-Mah-Tio}. 

\begin{figure}\label{excursion}
\begin{center}
\begin{tikzpicture}[scale=4]

\clip (1, 1) rectangle (-1, -0.3);

\draw (0, 0) circle (1);
\draw (0, 0.66666) circle (0.33333);

\filldraw[black] (0,0) circle (0.01);

\draw (0, 0) -- (0, 1);
\draw (0, 0) -- (0.5, 0.86603);
\draw[thick] (0, 0) -- (0.2, 0.97980);

\path (0,-0.1) node {$X_0$};
\path (-0.27,0.33) node {$H$};

\path (0.23,0.15) node {$\phi_{\text{max}}$};
\path (0.05,0.6) node {$\phi$};

\path (0.2,0.75) node {$\gamma$};
\path (-0.05,0.85) node {$\gamma_H$};
\path (0.5,0.7) node {$\gamma'_H$};

\draw[thick] (0,0) ++(60:0.2) arc (60:90:0.2);
\draw[thick] (0,0) ++(78:0.52) arc (78:90:0.52);

\end{tikzpicture}
\end{center} \caption{Excursion in the horoball $H$.}
\label{horoball}
\end{figure}
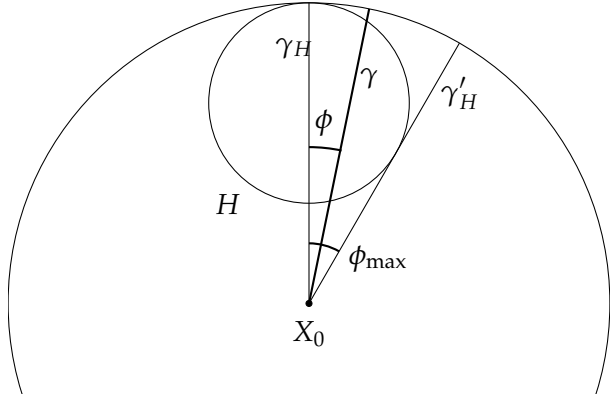

Let us fix the Poincare disc model for $\H(q)$.
In this model, the (unweighted by area) excursion $E(\gamma, H)$ can be interpreted as the \texttt{"}relative visual size\texttt{"} of the set of geodesic rays that go deeper than $\gamma$ in $H$. 
To be precise, suppose $\gamma$ is parameterised so that $\gamma(0) = X_0 \in \H(q)$. 
Let $\gamma_H$ be the geodesic ray from $X_0$ that goes straight to the cusp at infinity for the horoball $H$. 
Similarly, let $\gamma'_H$ be the geodesic ray from $X_0$ that is tangent to $H$ and on the same side of $\gamma_H$ as $\gamma$. 
Let $\phi_{\text{max}}$ be the angle between $\gamma_H$ and $\gamma'_H$ and $\phi$ be the angle between $\gamma_H$ and $\gamma$. 
See Figure \ref{excursion}.

By basic hyperbolic geometry, the excursion in $H$, unweighted by area, is up to uniform additive error given by 
\begin{equation}\label{e:excursion}
E(\gamma,H) = \frac{\phi_{\text{max}}}{\phi}.
\end{equation}
We leave the straightforward details of this to the reader.

Let $F^{\pm}$ be the vertical and horizontal measured foliations for $\gamma$. 
In the singular flat metric of the quadratic differential given by $\gamma_t$, let $\beta_t$ be a segment that is perpendicular to a core curve $\alpha$ of the cylinder. 
The \emph{twist parameter} $\text{tw}^-_t (\alpha)$ is the highest intersection number of a leaf of $F^-$ with $\beta_t$. 
Let $t_1$ and $t_2 \geqslant t_1$ be the times of entry and exit respectively of $\gamma$ in $H$.
The \emph{twists} $\text{tw}(\alpha)$ in $\alpha$ are defined as $\text{tw}^-_{t_2}(\alpha) - \text{tw}^-_{t_1}(\alpha) $. 
It does not matter that we use the horizontal foliation $F^-$. 
If we use the vertical foliation $F^+$, we get the same answer for $\text{tw}(\alpha)$ up to a uniform additive constant.  
At the end of the proof of \cite[Proposition 2.7]{Gad-Mah-Tio}, we prove that up to a uniform additive constant 
\begin{equation}\label{e:twist}
\text{tw}(\alpha) = \text{tw}^-_{t_2}(\alpha) - \text{tw}^-_{t_1}(\alpha) = \frac{2A}{\epsilon} \left(  \frac{\sin \phi_{\text{max}}} {\sin \phi}  \sqrt{1- \frac{\sin^2 \phi }{\sin^2 \phi_{\text{max}}}} \right)
\end{equation}
where $A$ is the flat area of the cylinder. 

\medskip
\noindent Let $\xi > 1$ and let $\xi'$ be the positive constant 
\[
\xi' = \frac{1}{2} \sqrt{1 - \frac{1}{\xi^2}}.
\]
Note that $\xi' < 1$. 
It is then elementary to verify that there is a constant $\phi_\xi> 0$ such that if $\phi_{\text{max}} < \phi_\xi$ and $\phi < \xi' \phi_{\text{max}}$ then
\[
\frac{\phi_{\text{max}}}{\xi \phi} < \frac{\sin \phi_{\text{max}}} {\sin \phi}  \sqrt{1- \frac{\sin^2 \phi }{\sin^2 \phi_{\text{max}}}} < \frac{\xi \phi_{\text{max}}}{\phi}.
\]
By relating $\phi_{\text{max}}$ to time along $\gamma'_H$ and hence $\gamma$, it follows that for any $\xi > 1$ there exists $s_\xi> 0$ such that if the entry time $t_1$ of $\gamma$ in $H$ satisfies $t_1 \geqslant s_\xi$ and if $\phi < \xi' \ \phi_{\text{max}}$ then 
\[
\frac{2A}{\epsilon \xi} \left( \frac{\phi_{\text{max}}}{\phi} \right) < \text{tw}(\alpha)  < \frac{2A \xi}{\epsilon} \left( \frac{\phi_{\text{max}}}{\phi} \right).
\]
By~(\ref{e:excursion}) the above bounds are equivalent to
\begin{equation}\label{e:unweighted}
\frac{2A}{\epsilon \xi} E(\gamma, H) < \text{tw}(\alpha) < \frac{2A\xi}{\epsilon} E(\gamma,H),
\end{equation}
provided $E(\gamma, H) > 1/\xi'$. 

\subsection{Proof of Theorem \ref{main}:} 
Let $\xi > 1$.
Consider~(\ref{e:unweighted}) for excursions that are larger than $1/\xi'$ and whose entry time is larger than $s_\xi$.
If we weight these excursions by area, that is consider $V_{\text{area}}$, then the relationship between twists and excursions becomes
\[
\frac{2}{\epsilon \xi} E_{V_{\text{area}}} (\gamma, H) < \text{tw}(\alpha) < \frac{2 \xi}{\epsilon} E_{V_{\text{area}}}(\gamma, H).
\]
We will use the relationship above in the excursion strong law~(\ref{t:strong}). 
Before we do that, we need to justify that the excursions we ignore do not affect the conclusion. 
This follows from the two observations below.
\begin{enumerate}
\item\label{bounded} Consider excursions $E(\gamma, H_j)$ that satisfy $E(\gamma, H_j) \leqslant 1/\xi'$ and consider their average as $T \to \infty$. 
By ergodicity of the \teichmuller geodesic flow, this average converges to the integral over $\mathcal{N}$ of a function that has finite expectation. 
See the discussion related to \cite[Inequality 3.6]{Gad}.
This implies that the sum of these shallow excursions is linear in $T$. 
Since the area of a cylinder is at most 1, we have $E_{V_{\text{area}}} (\gamma, H_j) \leqslant E(\gamma, H_j)$.
Hence the sum weighted by area of the shallow excursions is also linear in $T$.

\medskip
\item\label{finite} As $T \to \infty$, for all but finitely many initial excursions the entry times $t_j$ in the respective horoballs $H_j$ satisfy $t_j \geqslant s_\xi$. 
\end{enumerate}

By the observations above, we may ignore the excursions in~(\ref{bounded}) and~(\ref{finite}) in the strong law~(\ref{t:strong}).  
Thus, we may conclude that for any $\xi > 1$ there exists $T_\xi$ large enough depending on $\xi$ and the geodesic $\gamma$ such that 
\[
\left( \frac{2}{\epsilon \xi} \right) 2 \epsilon c_{\text{area}} (\mathcal{N}) <  \frac{\sum\limits_{j < N(T)} \text{tw}(\alpha_j)  - \max\limits_{j  < N(T)}\text{tw}(\alpha_j) }{T \log T} < \left( \frac{2 \xi}{ \epsilon }\right) 2 \epsilon c_{\text{area}}(\mathcal{N}) 
\]
for all $T > T_\xi$. As $\xi \to 1$, this implies
\[
\lim\limits_{T \to \infty} \frac{\sum\limits_{j < N(T)} \text{tw}(\alpha_j)  - \max\limits_{j < N(T) }\text{tw}(\alpha_j) }{T \log T} = 4 c_{\text{area}}
\]
finishing the proof of Theorem \ref{main}. 


\end{document}